\documentclass{article}
\usepackage{iclr2027_conference,times}
\usepackage[T1]{fontenc}
\usepackage{amsmath,amssymb,booktabs,graphicx,microtype,xcolor,array,tabularx}
\usepackage{subcaption}
\usepackage{tikz}
\usetikzlibrary{arrows.meta,positioning}
\usepackage{xurl}
\usepackage[breaklinks,colorlinks,citecolor=teal!70!black,linkcolor=teal!70!black,urlcolor=teal!70!black]{hyperref}
\newcommand{\FREQ}{\mathrm{FREQ}}
\newcommand{\WIE}{\mathrm{WIE}}
\newcommand{\NRE}{\mathrm{NRE}}
\newcommand{\fail}{\mathrm{X}}
\newcommand{\EERIE}{\texttt{EERIE}}
\title{\linespread{1.1}\selectfont Puzzle Game: Prediction and Classification of Wordle Solution Words}
\author{Haidong Xin\textsuperscript{*,\textdagger}\quad Fang Wu\textsuperscript{*}\quad Zhitong Zhou\textsuperscript{*}\\
Harbin Engineering University, China\\
\texttt{hdxin2002@gmail.com}}
\iclrfinalcopy
\begin{document}
\maketitle
\begingroup
\renewcommand{\thefootnote}{\fnsymbol{footnote}}
\footnotetext[1]{Equal contribution.}
\footnotetext[2]{Corresponding author.}
\endgroup
\lhead{Preprint}
\begin{abstract}
Wordle result prediction involves both daily participation and the difficulty of the solution word. We model these quantities with ARIMA forecasts of report counts, XGBoost regressors using word frequency, letter entropy, and repetition, and a decision tree trained on K-means difficulty tiers. Experiments use 355 daily outcome distributions, with a 107-word chronological holdout and four expanding-window tests. XGBoost reduces holdout mean absolute error from 4.08 to 3.84 percentage points and improves over the training mean in every temporal window. Ridge and random forests remain competitive. XGBoost systematically underestimates failure rates on difficult words. The lexical decision tree achieves 54.2\% tier accuracy, compared with 40.2\% when predicted histograms are assigned to the same training centroids. Bootstrap analysis shows that cluster stability does not imply that tiers are predictable from lexical attributes. For report counts, comparisons at matched horizons reveal substantial overdispersion in full-history ARIMA intervals. These results characterize the predictive value and limitations of a compact lexical model of reported Wordle performance.
\end{abstract}
\section{Introduction}
Wordle asks players to identify a five-letter word in at most six guesses, with feedback indicating correct letters and positions~\citep{nytwordle}. Publicly shared results record the proportions of reporting players who solve a puzzle in one through six guesses or fail. The daily number of reports reflects participation, whereas the outcome distribution reflects the difficulty of the solution word for that reporting population. Predicting these quantities requires accounting for both temporal variation and lexical attributes.

Predicting human performance differs from finding an effective Wordle strategy. Rank-one word representations~\citep{bonthron2022}, reinforcement learning and adaptive control~\citep{bhambri2022}, and information-based heuristics~\citep{greenberg2024} seek sequences of guesses that solve the puzzle efficiently. Human results average over different starting words, vocabularies, and strategies. Consequently, the performance of an algorithm with a fixed dictionary does not directly determine how difficult a word is for human players.

Previous studies have examined lexical neighborhoods and word difficulty~\citep{liu2023}, reporting dynamics and difficulty classification~\citep{weng2023}, and joint Bayesian prediction of reporting volume and attempt distributions~\citep{disilvio2023}. Behavioral analyses further show that first-guess success and public sharing depend on player behavior~\citep{dilger2023}. A useful prediction model must therefore specify its target population and evaluate performance on words that occur after its training period.

We consider three tasks: forecasting daily report counts, predicting outcome distributions from solution words, and assigning words to difficulty tiers (Figure~\ref{fig:framework}). ARIMA models reporting volume. XGBoost predicts the outcome histogram from word frequency (FREQ), within-word letter entropy (WIE), and the number of positions occupied by repeated letters (NRE). K-means groups training histograms into three tiers, and a shallow decision tree predicts those tiers from the same lexical attributes.

We evaluate the models on a chronological holdout and four expanding-window splits. The experiments measure the value of the lexical features, compare direct tier classification with classification through predicted histograms, and assess report-count forecasts at matched horizons. Error analysis identifies systematic underprediction of high failure rates. Resampling separates the stability of the difficulty clusters from the accuracy of lexical classification. An \EERIE{} case study illustrates the resulting distribution and tier predictions.

\begin{figure}[t]
\centering
\resizebox{\linewidth}{!}{%
\begin{tikzpicture}[
  box/.style={draw=black!50,rounded corners=2pt,minimum height=9mm,text width=25mm,align=center,font=\small,inner sep=3pt},
  method/.style={box,draw=blue!50!black,fill=blue!3},
  output/.style={box,draw=teal!70!black,fill=teal!3},
  flow/.style={-{Stealth[length=2mm]},line width=.6pt,draw=black!65},
  train/.style={flow,dashed}
]
\node[box] (counts) at (0,0) {Report Counts\\$r_t$};
\node[method] (arima) at (5.1,0) {ARIMA$(0,1,1)$};
\node[output] (forecast) at (10.2,0) {Count Forecast\\$\widehat r_{T+h\mid T}$};
\draw[flow] (counts) -- (arima);
\draw[flow] (arima) -- (forecast);
\node[box] (word) at (0,-1.65) {Solution Word\\$w$};
\node[method] (features) at (3.4,-1.65) {Lexical Attributes\\$x(w)$};
\node[method] (boost) at (6.8,-1.65) {XGBoost\\Normalization};
\node[output] (hist) at (10.2,-1.65) {Histogram\\$\widehat{\mathbf p}(w)$};
\draw[flow] (word) -- (features);
\draw[flow] (features) -- (boost);
\draw[flow] (boost) -- (hist);
\node[box] (training) at (0,-3.3) {Training Data\\$\mathbf p_i,\;i\in\mathcal T$};
\node[method] (cluster) at (3.4,-3.3) {K-means\\$K=3$};
\node[method] (tree) at (6.8,-3.3) {Decision Tree\\$g(x)$};
\node[output] (tier) at (10.2,-3.3) {Difficulty Tier\\$g(x(w))$};
\draw[train] (training) -- (cluster);
\draw[train] (cluster) -- node[above,font=\footnotesize] {$c_i$} (tree);
\draw[flow] (features.south) -- (tree.north);
\draw[flow] (tree) -- (tier);
\end{tikzpicture}%
}
\caption{Reporting-volume and word-level prediction. Dashed arrows construct the training labels for the decision tree; solid arrows describe prediction inputs and outputs.}
\label{fig:framework}
\end{figure}

\section{Data Preparation}
\label{sec:data}
\paragraph{Data Source.}
The 2023 MCM Problem C dataset contains 359 daily records from January 7 to December 31, 2022~\citep{comap2023}. Each record includes a date, puzzle index, solution word, total number of reported results, hard-mode report count, and seven outcome percentages. The records were collected from results shared on Twitter. Our prediction target is therefore the outcome distribution among \emph{reporting players}.

\paragraph{Cleaning.}
We trim surrounding whitespace, uppercase words, and remove diacritics, so that \texttt{favor } and \texttt{na\"ive} become \texttt{FAVOR} and \texttt{NAIVE}. The entries \texttt{TASH}, \texttt{CLEN}, and \texttt{RPROBE} fail the five-letter schema and are excluded from word-level modeling. The March 27 histogram sums to 126\% and is also excluded. For seven integer-rounded shares, the largest possible total rounding discrepancy is 3.5 percentage points (pp). We retain rows within that tolerance and normalize each retained vector to 100\%. The resulting word-level dataset contains 355 observations.

The count series retains all 359 dates, because an invalid word string does not invalidate a daily report count. On November 30, the \emph{total} report count is 2,569; the hard-mode count is 2,405. We replace this isolated anomaly with the mean of the preceding six daily totals, 25,546.5. Appendix~\ref{app:forecasts} reports sensitivity to using the recorded total instead.

\paragraph{Prediction Targets.}
Let $\mathcal{B}=\{1,2,3,4,5,6,\fail\}$ index the seven outcomes. For word $w_i$, the vector $\mathbf{p}_i=(p_{i,b})_{b\in\mathcal{B}}$ contains percentages satisfying $p_{i,b}\geq0$ and $\sum_{b\in\mathcal{B}}p_{i,b}=100$. The category $\fail$ denotes failure within six guesses. The calendar series $r_t$ is the number of reports on day $t$. Word-level observations receive equal fitting weight, regardless of their daily report counts.

\paragraph{Evaluation Protocol.}
The index sets $\mathcal{T}$ and $\mathcal{E}$ contain the first 248 and last 107 retained words, respectively. The training partition ends on September 13, 2022; the test partition starts on September 14. Hyperparameters are selected using the last 50 training words as validation data and the preceding 198 for fitting. Each selected model is then refitted on all of $\mathcal{T}$. Test outcomes enter neither model selection nor cluster fitting. A separate count-series evaluation uses November--December as a 61-day sequential holdout. These chronological protocols match the future-prediction setting; random-split sensitivity answers a different generalization question~\citep{bergmeir2018,hyndman2021}.

We also evaluate four expanding-window splits, with training sizes 155, 205, 255, and 305 and a 50-word test block following each training period. Every split repeats selection within its own training prefix, using its last 20\% for validation. The test blocks are disjoint; training sets overlap. The primary holdout and expanding-window tests evaluate temporal generalization within 2022~\citep{roberts2017,cawley2010}. Appendix~\ref{app:protocol} gives exact dates and candidate settings.

\section{Methodology}
\subsection{Reporting Volume Forecasting}
We model report counts with the ARIMA$(0,1,1)$ specification~\citep{newbold1983,hyndman2021}
\begin{equation}
 (1-B)r_t=\varepsilon_t+\theta\varepsilon_{t-1},
 \qquad \varepsilon_t\overset{\mathrm{iid}}{\sim}\mathcal{N}(0,\sigma^2),
 \label{eq:arima}
\end{equation}
where $B$ is the backshift operator defined by $Br_t=r_{t-1}$, $\theta\in(-1,1)$ is the moving-average coefficient, and $\sigma^2>0$ is the innovation variance. The model has no drift term. We estimate its parameters by state-space maximum likelihood in statsmodels~\citep{seabold2010}, using the same sign convention for $\theta$ in estimation and sensitivity analysis.

For sequential evaluation, parameters are fitted through October 31. Each subsequent one-step forecast is issued before observing that day's count; the observation then updates the model state without re-estimating parameters. The baseline predicts the previous day's count. We report MAE, MAPE, and empirical coverage of model-based 95\% prediction intervals.

The March 1, 2023 forecast has origin $T=$ December 31, 2022 and horizon $h=60$. Under the fitted Gaussian model, its interval is
\begin{equation}
 I_{T,h}^{1-\alpha}
 =\bigl[m_{T,h}-z_{1-\alpha/2}\sqrt{v_{T,h}},\;
         m_{T,h}+z_{1-\alpha/2}\sqrt{v_{T,h}}\bigr],
 \label{eq:interval}
\end{equation}
where $m_{T,h}=\widehat r_{T+h\mid T}$ and $v_{T,h}$ are the conditional forecast mean and variance given counts through $T$, with fitted parameters held fixed, and $z_{1-\alpha/2}$ is a standard-normal quantile. We set $\alpha=0.05$. The interval incorporates multi-step forecast variance but excludes parameter-estimation uncertainty. We retain negative lower endpoints to expose the Gaussian count-model limitation. Augmented Dickey--Fuller~\citep{dickey1979} and Ljung--Box~\citep{ljung1978} tests provide additional diagnostics.

\subsection{Lexical Feature Representation}
\label{sec:features}
We use three lexical attributes: word frequency, within-word letter entropy, and repetition. Corpus frequency is a familiarity proxy whose interpretation depends on register~\citep{brysbaert2009,michel2011}. FREQ is the English probability returned by \texttt{wordfreq} 3.1.1~\citep{wordfreq}. We use the same fixed English frequency estimates for every observation date.

For a word $w=w_1\cdots w_5$, define $n_c(w)=\sum_{j=1}^{5}\mathbf{1}\{w_j=c\}$ and $q_c(w)=n_c(w)/5$. The two structural attributes are
\begin{equation}
 \begin{aligned}
 \WIE(w)&=-\sum_{c\in\mathcal{A}:q_c(w)>0}q_c(w)\log_2q_c(w),\\
 \NRE(w)&=\sum_{c\in\mathcal{A}}n_c(w)\mathbf{1}\{n_c(w)\geq2\},
 \end{aligned}
 \label{eq:features}
\end{equation}
where $\mathcal{A}$ is the English alphabet, $n_c(w)$ counts occurrences of letter $c$, $q_c(w)$ is its within-word proportion, and $\mathbf{1}\{\cdot\}$ is an indicator. WIE is measured in bits and lies in $[0,\log_2 5]$~\citep{shannon1948}. NRE counts positions occupied by repeated letters, rather than repeated letter types or excess occurrences. For example, \EERIE{} has $\WIE=1.37095$ and $\NRE=3$.

The input representation is
\begin{equation}
 x(w)=\bigl(\log_{10}\max\{\FREQ(w),\epsilon\},\;\WIE(w),\;\NRE(w)\bigr)^{\!\top}\in\mathbb{R}^{3},
 \label{eq:representation}
\end{equation}
where $\epsilon=10^{-12}$ prevents taking the logarithm of zero. The floored log transform is monotone in FREQ. WIE describes letter diversity; it is not the expected information gain of a guess, which depends on a candidate set and feedback patterns~\citep{disilvio2023,greenberg2024}. WIE and NRE encode related repetition patterns, so their individual importances need not represent independent effects.

In the descriptive profiles, high FREQ denotes a value at or above the dataset median, high WIE denotes five distinct letters, and high NRE denotes any repeated letter (Figure~\ref{fig:profiles}, Appendix~\ref{app:errors}). Feature ablations assess these attribute groups under fixed model capacity.

\subsection{Outcome Distribution Prediction}
\label{sec:xgb}
For each $b\in\mathcal{B}_{+}=\{2,3,4,5,6,\fail\}$, we fit an independent XGBoost regressor~\citep{friedman2001,chen2016}. Its stage-wise update is
\begin{equation}
 z_b^{(0)}(x)=\overline p_b,\qquad
 z_b^{(m)}(x)=z_b^{(m-1)}(x)+\eta f_{m,b}(x),
 \label{eq:boosting}
\end{equation}
where $\mathcal I$ is the current fitting index set, $\overline p_b=|\mathcal I|^{-1}\sum_{i\in\mathcal I}p_{i,b}$ is its outcome mean, $\eta\in(0,1]$ is the learning rate, $m\in\{1,\ldots,M\}$ indexes boosting rounds, and $f_{m,b}$ is a regression tree. Before shrinkage by $\eta$, each added tree approximately minimizes
\begin{equation}
 \mathcal{L}_{b}^{(m)}(f)
 =\frac12\sum_{i\in\mathcal I}
   \bigl(p_{i,b}-z_b^{(m-1)}(x_i)-f(x_i)\bigr)^2
 +\gamma L(f)+\frac{\lambda}{2}\|a(f)\|_2^2,
 \label{eq:xgb-objective}
\end{equation}
where $x_i=x(w_i)$, $L(f)$ is the tree's number of leaves, $a(f)$ is its vector of leaf weights, and $\gamma,\lambda\geq0$ penalize tree complexity and leaf magnitudes. 
After $M$ rounds, we use $z_b=z_b^{(M)}$ for the six learned bins. The low-mass one-guess outcome uses $z_1(x)=\overline p_1$.

Separate regressors need not produce a valid composition~\citep{aitchison1982}. Let $u_b(x)=[z_b(x)]_+$ and $S(x)=\sum_{j\in\mathcal{B}}u_j(x)$. The final percentages are
\begin{equation}
 \widehat p_b(x)=
 \begin{cases}
 100\,u_b(x)/S(x),&S(x)>0,\\
 100/|\mathcal{B}|,&S(x)=0,
 \end{cases}
 \qquad b\in\mathcal{B},
 \label{eq:normalization}
\end{equation}
where $[z]_+=\max\{z,0\}$ is the positive part and $S(x)$ is the total clipped mass. Thus $\widehat p_b(x)\geq0$ and $\sum_{b\in\mathcal{B}}\widehat p_b(x)=100$. Normalization can change the initially constant first-bin share.

We select depth from $\{1,2,3\}$ and rounds from $\{100,300\}$ using validation MAE after normalization. Fixed settings are $\eta=0.03$, $\lambda=10$, $\gamma=0$, minimum child weight 5, and full row and column sampling. The fitting set contains 198 words during primary selection and 248 for the evaluated model; the selected model uses depth 1 and 300 rounds. Comparators are the training-mean histogram, standardized ridge regression~\citep{hoerl1970}, distance-weighted nearest-neighbor regression~\citep{stone1977}, and a multi-output random forest~\citep{breiman2001}. All use the same training/validation boundaries and three lexical inputs; learned predictions use Equation~\eqref{eq:normalization}. Appendix~\ref{app:protocol} specifies their grids. All stochastic procedures use seed 42 unless otherwise stated.

\subsection{Difficulty Tier Classification}
\label{sec:tiers}
We form three difficulty tiers by applying K-means to training outcome percentages~\citep{lloyd1982}. The clustering objective is
\begin{equation}
 \min_{\mu_1,\ldots,\mu_K\in\mathbb{R}^{7}}
 \sum_{i\in\mathcal{T}}\min_{1\leq k\leq K}
 \|\mathbf p_i-\mu_k\|_2^2,
 \qquad K=3,
 \label{eq:clustering}
\end{equation}
where $\mu_k$ is a cluster centroid and Euclidean distance is applied without per-bin variance scaling. We use k-means++ initialization~\citep{arthur2007} and 50 restarts. The centroids are relabeled in increasing order of the capped score $d(\mu)=\bigl(\sum_{b=1}^{6}b\mu_b+7\mu_{\fail}\bigr)/100$ and named easy, moderate, and difficult. The division by 100 converts percentages to probabilities. Coding failure as 7 defines an ordering convention; it does not estimate failed players' eventual guess counts.

After relabeling, an observed histogram receives target tier $c_i=\arg\min_{k\in\{1,2,3\}}\|\mathbf p_i-\mu_k\|_2^2$. Ties, if any, follow a fixed implementation order. A decision tree $g(x_i)$ predicts this tier from lexical attributes~\citep{breiman1984}, with maximum depth 3 and at least 10 training words per leaf. Test outcomes determine evaluation labels only, using the fixed training centroids. We report accuracy, balanced accuracy, macro-F1, and a training-majority baseline. The tiers summarize empirical outcome profiles. We examine their separation using elbow and silhouette diagnostics~\citep{rousseeuw1987}, and their stability by resampling the training profiles.

\section{Experimental Results}
\subsection{Distribution Prediction Under Temporal Shift}
Table~\ref{tab:prediction} compares five predictors on the primary 107-word holdout. XGBoost reduces MAE by 0.243 pp, or 6.0\%, relative to the training mean. Its paired seven-word block-bootstrap interval for the difference is $[-0.424,-0.036]$ pp. Ridge and random forests improve on the mean as well; the forest attains the lowest point errors. The forest-minus-XGBoost MAE difference is $-0.089$ pp with conditional interval $[-0.175,-0.009]$. 

\begin{table}[t]
\caption{Primary chronological holdout ($n=107$). JSD uses bits; $W_1$ uses the capped outcome order $1,\ldots,6,\fail$. Intervals are paired circular-block bootstrap intervals for MAE differences in pp, with block length 7 and 2,000 draws. Best point errors are bold.}
\label{tab:prediction}\centering\small
\begin{tabular*}{\linewidth}{@{\extracolsep{\fill}}lrrrr@{}}
\toprule
Model & MAE (pp) $\downarrow$ & JSD $\downarrow$ & $W_1$ $\downarrow$ & $\Delta$MAE vs. Mean (95\% CI)\\
\midrule
Training Mean & 4.082 & 0.0330 & 0.366 & --\\
Ridge & 3.799 & 0.0294 & 0.333 & $[-0.490,\,-0.062]$\\
kNN & 4.262 & 0.0363 & 0.380 & $[-0.120,\,+0.490]$\\
Random Forest & \textbf{3.750} & \textbf{0.0288} & \textbf{0.331} & $[-0.522,\,-0.134]$\\
XGBoost & 3.839 & 0.0297 & 0.337 & $[-0.424,\,-0.036]$\\
\bottomrule
\end{tabular*}

\end{table}

MAE measures bin-wise percentage errors, JSD compares probability vectors~\citep{lin1991}, and $W_1$ measures displacement along the capped outcome order~\citep{villani2009}. TV is exactly $7\operatorname{MAE}_{\mathrm{pp}}/200$ under our weighting, so it supplies no independent evidence. Bin-wise $\pm3$ pp accuracy is reported separately in Appendix~\ref{app:metrics}; it is not simultaneous coverage of a new histogram.

\begin{table}[t]
\caption{Nested expanding-window tests. Each fold predicts the next 50 words after training on 155, 205, 255, or 305 words. Entries are MAE (pp); the final column is the mean and sample SD across four folds, not a confidence interval.}
\label{tab:rolling}\centering\small
\begin{tabular*}{\linewidth}{@{\extracolsep{\fill}}lrrrrr@{}}
\toprule
Model & F1 & F2 & F3 & F4 & Mean $\pm$ SD\\
\midrule
Training Mean & 3.610 & 4.279 & 4.077 & 3.724 & $3.923\pm0.310$\\
Ridge & 3.084 & 3.808 & 3.788 & 3.560 & $3.560\pm0.337$\\
kNN & 3.605 & 4.148 & 4.136 & 4.106 & $3.999\pm0.263$\\
Random Forest & 3.160 & 3.884 & 3.651 & 3.605 & $3.575\pm0.302$\\
XGBoost & 3.176 & 3.887 & 3.883 & 3.547 & $3.623\pm0.338$\\
\bottomrule
\end{tabular*}

\end{table}

\begin{figure}[t]
\centering
\begin{subfigure}[t]{.32\linewidth}\centering
\includegraphics[width=\linewidth]{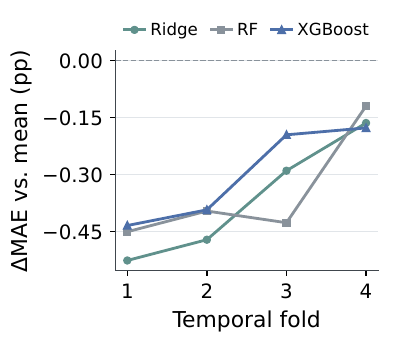}
\caption{Temporal improvements}
\end{subfigure}\hfill
\begin{subfigure}[t]{.32\linewidth}\centering
\includegraphics[width=\linewidth]{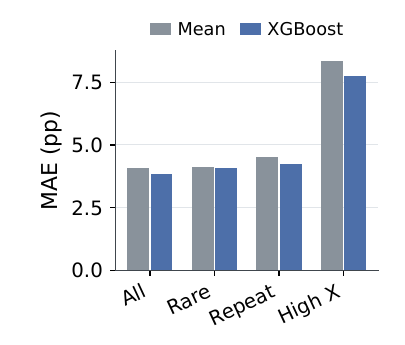}
\caption{Error concentration}
\end{subfigure}\hfill
\begin{subfigure}[t]{.32\linewidth}\centering
\includegraphics[width=\linewidth]{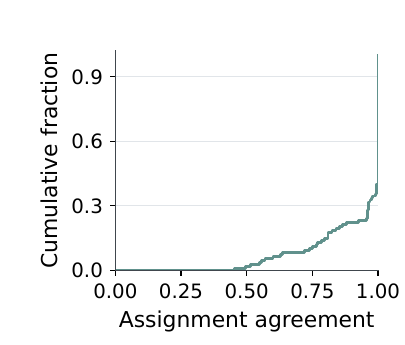}
\caption{Tier stability}
\end{subfigure}
\caption{Temporal and structural robustness. Negative differences in (a) favor the learned model. Subsets in (b) contain 107, 55, 26, and 8 words. Panel (c) shows per-word agreement with the reference tier over 200 training-profile bootstrap resamples.}
\label{fig:robustness}
\end{figure}

Table~\ref{tab:rolling} tests whether the primary result depends on one cutoff. XGBoost improves over the training mean in all four windows (Figure~\ref{fig:robustness}a), with average MAE 3.623 versus 3.923 pp. Ridge and random forests average 3.560 and 3.575 pp. Distance-weighted kNN averages 3.999 pp, suggesting that proximity in this compact feature space is insufficient on its own. Fold-to-fold variation is larger than the mean difference between ridge and boosting. Because training prefixes overlap, we do not treat the fold SD as a standard error or conduct an independent-fold significance test~\citep{bengio2004}. The improvements over the mean are smaller in the later windows.

\subsection{Error Analysis}
The largest errors concern unusual outcome shapes. On 26 test words with repeated letters, XGBoost MAE is 4.222 pp, compared with 3.716 pp on 81 all-distinct words. For the eight words whose observed failure share is at or above the training 90th-percentile threshold of 6.06\%, MAE rises to 7.735 pp and the mean failure-bin bias is $-15.25$ pp (Figure~\ref{fig:robustness}b). The training mean also struggles on this subset (8.333 pp). Negative tail bias persists in each expanding window under its own training-derived threshold. The reduction in average error therefore coexists with substantial underprediction on difficult words. Appendix~\ref{app:errors} reports the window-specific counts, all model slices, signed errors, and largest-error cases.

Repetition-related attributes carry much of the usable structure. Fixed-capacity XGBoost ablations yield MAE 4.124 pp with FREQ alone, 3.821 with WIE+NRE, and 3.839 with all three. Jointly permuting WIE and NRE in the fixed test inputs increases MAE by 0.306 pp on average across 200 permutations; permuting FREQ increases it by 0.108 pp. Joint permutation preserves the internal relationship of the two repetition measures, but may still create unrealistic combinations with frequency. These permutation effects measure reliance on the fitted features~\citep{strobl2008}. Under the selected model capacity, adding FREQ does not improve on WIE+NRE in the chronological test.

\subsection{Difficulty Classification And Stability}
\begin{figure}[t]
\centering
\begin{subfigure}[t]{0.32\linewidth}
\centering
\includegraphics[width=\linewidth]{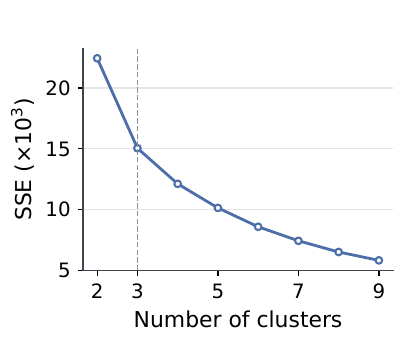}
\caption{Cluster distortion}\label{fig:tiers:a}
\end{subfigure}
\hfill
\begin{subfigure}[t]{0.32\linewidth}
\centering
\includegraphics[width=\linewidth]{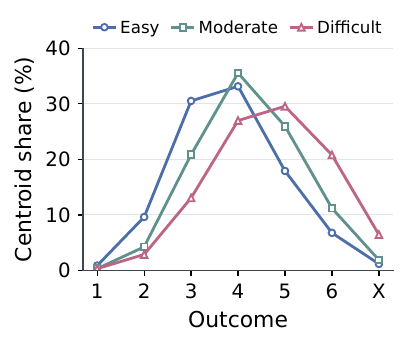}
\caption{Cluster centroids}\label{fig:tiers:b}
\end{subfigure}
\hfill
\begin{subfigure}[t]{0.32\linewidth}
\centering
\includegraphics[width=\linewidth]{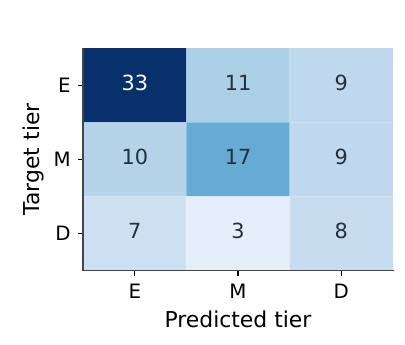}
\caption{Tier classification}\label{fig:tiers:c}
\end{subfigure}
\caption{Difficulty clusters and decision-tree predictions on held-out words. E, M, and D denote easy, moderate, and difficult.}
\label{fig:tiers}
\end{figure}

The training clusters contain 77 easy, 107 moderate, and 64 difficult words. Their centroids shift mass toward later outcomes (Figure~\ref{fig:tiers}b). The silhouette score is 0.371 for $K=3$, versus 0.418 for $K=2$; the three-tier vocabulary is an interpretive choice rather than evidence of three uniquely natural populations.

We refit K-means on 200 bootstrap resamples of the 248 training profiles and reassign the same 107 test profiles. The mean adjusted Rand index against the reference partition is 0.816, with a 2.5th--97.5th percentile resampling range of $[0.513,0.980]$~\citep{hubert1985}. Ordered-label agreement averages 93.0\%; 16 words agree with their reference tier in fewer than 80\% of resamples (Figure~\ref{fig:robustness}c). Most profiles are reproducibly assigned, with instability concentrated near centroid boundaries (Appendix~\ref{app:stability}).

The lexical decision tree attains 54.2\% accuracy, 51.3\% balanced accuracy, and 0.504 macro-F1, versus 33.6\% majority accuracy. Assigning XGBoost's predicted histograms to the same training centroids instead gives 40.2\% accuracy: 90 words are predicted moderate, although only 36 reference labels are moderate. A low histogram loss therefore need not yield good tier decisions. The route comparison and a sufficient distance-margin condition are detailed in Appendix~\ref{app:routes}. Despite their resampling stability, the outcome clusters remain only partly predictable from these features. Sensitivity to $K\in\{2,3,4,5\}$ is reported without choosing $K$ by test performance.

\subsection{Report Forecasts Across Matched Horizons}
November--December one-step forecasts retain parameters fitted through October 31 and update states sequentially. ARIMA MAE is 1,547 reports versus 1,743 for persistence; MAPE is 6.55\% versus 7.37\%. Its 100\% empirical 95\% interval coverage masks excessive dispersion.

For horizon-matched evaluation, 18 weekly origins from July 4 to October 31 each evaluate $h\in\{1,7,14,30,60\}$. Fixed sensitivity models apply the same ARIMA order to full history, trailing 90-day counts, or log counts within that window, separating estimation-window and support effects.

\begin{figure}[t]
\centering
\begin{subfigure}[t]{.32\linewidth}\centering
\includegraphics[width=\linewidth]{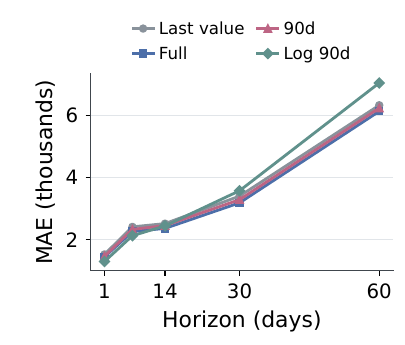}
\caption{Point errors}
\end{subfigure}\hfill
\begin{subfigure}[t]{.32\linewidth}\centering
\includegraphics[width=\linewidth]{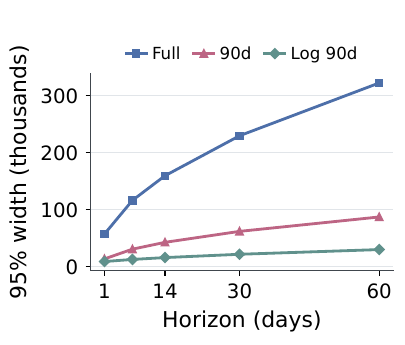}
\caption{Interval width}
\end{subfigure}\hfill
\begin{subfigure}[t]{.32\linewidth}\centering
\includegraphics[width=\linewidth]{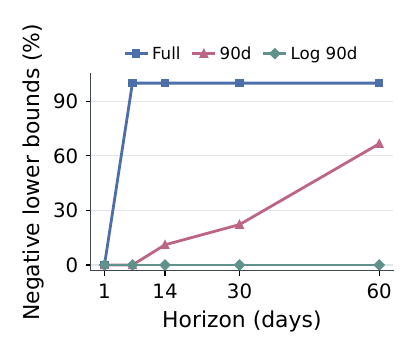}
\caption{Count-support violations}
\end{subfigure}
\caption{Horizon-matched report forecasts at 18 weekly origins. Full denotes full-history ARIMA; 90d and Log 90d use the trailing 90 days. All ARIMA models have order $(0,1,1)$ with no drift. Log-scale means include the Gaussian variance correction.}
\label{fig:horizons}
\end{figure}

At $h=60$, full-history ARIMA has MAE 6,128 reports versus 6,323 for persistence, but mean 95\% interval width is 322,639 reports. The 90-day window reduces width to 87,168; log-count ARIMA gives positive endpoints and width 29,880, yet its bias-corrected mean has larger MAE, 7,040. Positivity and sharpness therefore do not ensure better point prediction. Figure~\ref{fig:horizons} and Appendix~\ref{app:forecasts} report horizon-specific point and interval scores~\citep{gneiting2007}. Overlapping origins in one decline period cannot establish calibration on a future year.

Fitting ARIMA to all 359 cleaned counts gives a March 1, 2023 forecast of 20,521 with interval $[-109{,}546,150{,}588]$. The negative endpoint exposes misspecification; further diagnostics and window sensitivities appear in Appendix~\ref{app:diagnostics}.

\subsection{EERIE Case Study}
\begin{figure}[!ht]
\centering
\begin{subfigure}[t]{0.32\linewidth}
\centering
\includegraphics[width=\linewidth]{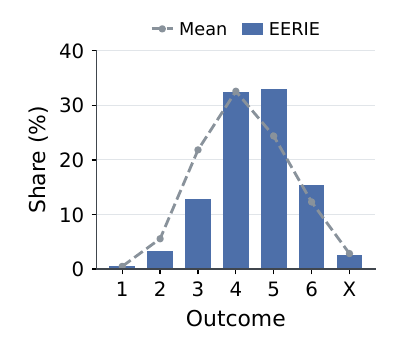}
\caption{EERIE prediction}\label{fig:eerie:a}
\end{subfigure}
\hfill
\begin{subfigure}[t]{0.32\linewidth}
\centering
\includegraphics[width=\linewidth]{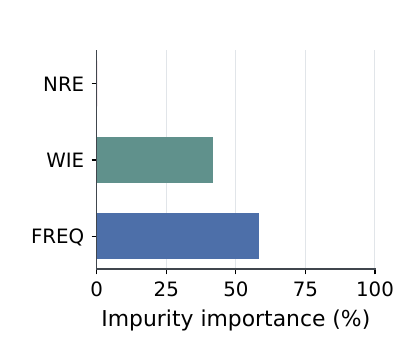}
\caption{Feature importance}\label{fig:eerie:b}
\end{subfigure}
\hfill
\begin{subfigure}[t]{0.32\linewidth}
\centering
\includegraphics[width=\linewidth]{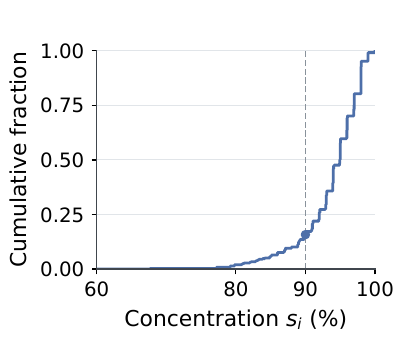}
\caption{Outcome concentration}\label{fig:eerie:c}
\end{subfigure}
\caption{Predicted outcomes for \EERIE{} from the 248-word model, decision-tree feature importance, and the empirical cumulative distribution of outcome concentration. The vertical line in (c) marks 90\%.}
\label{fig:eerie}
\end{figure}

The MCM query \EERIE{}~\citep{comap2023} has $\FREQ=2.14\times10^{-6}$, $\WIE=1.37095$, and $\NRE=3$. The primary 248-word XGBoost model predicts $(0.48,3.34,12.85,32.42,32.93,15.34,2.64)$\% for outcomes $(1,2,3,4,5,6,\fail)$. Both its histogram and the lexical tree assign the difficult tier. Their shared features and tier definition limit what can be inferred from this agreement; the case has no observed evaluation target. Figure~\ref{fig:eerie} and Appendix~\ref{app:cases} give the case predictions, including the all-data refit, without an additional test-performance claim.

\section{Discussion}
The three lexical attributes improve distribution prediction over a training-mean histogram in all four temporal windows. Similar or lower errors from ridge and random forests indicate that the predictive signal is not specific to boosting. Repetition accounts for much of this signal, but the features do not distinguish letter identities or positions. Words with similar feature values can therefore have substantially different outcome distributions.

The dataset contains 255 words with distinct letters, 94 with one pair, four with two pairs, and two with a triple. On these patterns, WIE and NRE determine one another. Their different numerical orderings can affect tree partitions, but using both does not add distinctions between words. The representation therefore remains too coarse to describe variation within each repetition pattern.

Two findings limit the interpretation of average histogram error. First, near-zero aggregate failure bias conceals underprediction on high-failure words in every temporal window. Second, assigning predicted histograms to the training centroids concentrates predictions in the moderate tier. Direct classification gives higher tier accuracy, suggesting that the choice of prediction target matters even when the input features are unchanged.

The report-count results show a separate modeling limitation. Early participation volatility produces wide long-horizon intervals under full-history ARIMA. A recent estimation window reduces their width, while a log transform enforces positive forecasts but can increase point error. Forecast assessment should therefore consider point loss, interval coverage, and interval width together.

\paragraph{Limitations.}
The observations cover one year and a self-selected reporting population. Starting guesses, player vocabularies, and lexical neighborhoods are unobserved, and the fixed frequency estimates may differ from usage in 2022. The temporal evaluations share training observations, and bootstrap intervals condition on fitted models. Evaluation on later years and richer word representations would help establish how well the results transfer.

\section{Conclusion}
We studied Wordle reporting volume, outcome distributions, and difficulty tiers using ARIMA, XGBoost, and K-means with a lexical decision tree. Chronological tests show useful lexical signal alongside systematic errors on high-failure words. Direct tier classification improves on classification through predicted histograms, while report-count forecasts expose the sensitivity of ARIMA intervals to the estimation window. These findings motivate models that incorporate lexical neighborhoods and player strategies, with evaluation across later puzzle periods.

\clearpage
\begingroup
\small
\bibliographystyle{iclr2027_conference}
\bibliography{references}

@misc{comap2023,
 author={{COMAP}}, title={2023 {MCM} Problem {C}: Predicting {Wordle} Results}, year={2023},
 howpublished={Mathematical Contest in Modeling problem statement and accompanying data},url={https://www.comap.org/component/zoo/?element=71992483-3d3b-47dc-a146-3c901bdea04f&format=raw&item_id=3324&method=download&task=callelement}}

@misc{nytwordle,author={{The New York Times}},title={{Wordle}},year={2022},url={https://www.nytimes.com/games/wordle/index.html}}

@misc{bonthron2022,author={Bonthron, Michael},title={Rank One Approximation as a Strategy for {Wordle}},year={2022},note={arXiv:2204.06324},url={https://arxiv.org/abs/2204.06324}}

@misc{bhambri2022,author={Bhambri, Siddhant and Bhattacharjee, Amrita and Bertsekas, Dimitri},title={Reinforcement Learning Methods for {Wordle}: A {POMDP}/Adaptive Control Approach},year={2022},note={arXiv:2211.10298},url={https://arxiv.org/abs/2211.10298}}

@misc{greenberg2024,author={Greenberg, Ronald I.},title={Effective {Wordle} Heuristics},year={2024},note={arXiv:2408.11730, version 2, revised 2026},url={https://arxiv.org/abs/2408.11730v2}}

@misc{liu2023,author={Liu, Beibei and Zhang, Yuanfang and Zhang, Shiyu},title={Explore the Difficulty of Words and Its Influential Attributes Based on the {Wordle} Game},year={2023},note={arXiv:2305.03502},url={https://arxiv.org/abs/2305.03502}}

@misc{weng2023,author={Weng, Jiaqi and Feng, Chunlin},title={Prediction Model for {Wordle} Game Results with High Robustness},year={2023},note={arXiv:2309.14250},url={https://arxiv.org/abs/2309.14250}}

@misc{disilvio2023,author={DiSilvio, Steven and Ozerov, Anthony and Zhou, Leon},title={How Many {Wordle} Words Will {Wordle} Guessers Guess If {Wordle}'s Wednesday {Wordle} Word Is ``Eerie''?},year={2023},note={arXiv:2311.16777},url={https://arxiv.org/abs/2311.16777}}

@misc{dilger2023,author={Dilger, James P.},title={{Wordle}: A Microcosm of Life. Luck, Skill, Cheating, Loyalty, and Influence!},year={2023},note={arXiv:2309.02110},url={https://arxiv.org/abs/2309.02110}}

@article{shannon1948,author={Shannon, Claude E.},title={A Mathematical Theory of Communication},journal={Bell System Technical Journal},volume={27},number={3},pages={379--423},year={1948},url={https://doi.org/10.1002/j.1538-7305.1948.tb01338.x}}

@article{brysbaert2009,author={Brysbaert, Marc and New, Boris},title={Moving beyond {Ku\v{c}era} and {Francis}: A Critical Evaluation of Current Word Frequency Norms and the Introduction of a New and Improved Word Frequency Measure for {American English}},journal={Behavior Research Methods},volume={41},number={4},pages={977--990},year={2009},url={https://doi.org/10.3758/BRM.41.4.977}}

@article{michel2011,author={Michel, Jean-Baptiste and Shen, Yuan Kui and Aiden, Aviva Presser and Veres, Adrian and Gray, Matthew K. and {The Google Books Team} and Pickett, Joseph P. and Hoiberg, Dale and Clancy, Dan and Norvig, Peter and Orwant, Jon and Pinker, Steven and Nowak, Martin A. and Aiden, Erez Lieberman},title={Quantitative Analysis of Culture Using Millions of Digitized Books},journal={Science},volume={331},number={6014},pages={176--182},year={2011},url={https://doi.org/10.1126/science.1199644}}

@misc{wordfreq,author={Speer, Robyn},title={rspeer/wordfreq: v3.0},year={2022},publisher={Zenodo},note={Recommended software citation, v3.0.2; experiments use v3.1.1},url={https://doi.org/10.5281/zenodo.7199437}}

@article{newbold1983,author={Newbold, Paul},title={{ARIMA} Model Building and the Time Series Analysis Approach to Forecasting},journal={Journal of Forecasting},volume={2},number={1},pages={23--35},year={1983},url={https://doi.org/10.1002/for.3980020104}}

@book{hyndman2021,author={Hyndman, Rob J. and Athanasopoulos, George},title={Forecasting: Principles and Practice},edition={3},publisher={OTexts},address={Melbourne, Australia},year={2021},url={https://otexts.com/fpp3/}}

@article{dickey1979,author={Dickey, David A. and Fuller, Wayne A.},title={Distribution of the Estimators for Autoregressive Time Series with a Unit Root},journal={Journal of the American Statistical Association},volume={74},number={366},pages={427--431},year={1979},url={https://doi.org/10.2307/2286348}}

@article{ljung1978,author={Ljung, Greta M. and Box, George E. P.},title={On a Measure of Lack of Fit in Time Series Models},journal={Biometrika},volume={65},number={2},pages={297--303},year={1978},url={https://doi.org/10.1093/biomet/65.2.297}}

@article{friedman2001,author={Friedman, Jerome H.},title={Greedy Function Approximation: A Gradient Boosting Machine},journal={The Annals of Statistics},volume={29},number={5},pages={1189--1232},year={2001},url={https://doi.org/10.1214/aos/1013203451}}

@inproceedings{chen2016,author={Chen, Tianqi and Guestrin, Carlos},title={{XGBoost}: A Scalable Tree Boosting System},booktitle={Proceedings of the 22nd ACM SIGKDD International Conference on Knowledge Discovery and Data Mining},pages={785--794},year={2016},url={https://doi.org/10.1145/2939672.2939785}}

@article{lloyd1982,author={Lloyd, Stuart P.},title={Least Squares Quantization in {PCM}},journal={IEEE Transactions on Information Theory},volume={28},number={2},pages={129--137},year={1982},url={https://doi.org/10.1109/TIT.1982.1056489}}

@inproceedings{arthur2007,author={Arthur, David and Vassilvitskii, Sergei},title={k-means++: The Advantages of Careful Seeding},booktitle={Proceedings of the Eighteenth Annual ACM-SIAM Symposium on Discrete Algorithms},pages={1027--1035},year={2007},url={https://theory.stanford.edu/~sergei/papers/kMeansPP-soda.pdf}}

@book{breiman1984,author={Breiman, Leo and Friedman, Jerome H. and Olshen, Richard A. and Stone, Charles J.},title={Classification and Regression Trees},publisher={Wadsworth},address={Belmont, CA},year={1984},url={https://www.routledge.com/Classification-and-Regression-Trees/Breiman-Friedman-Olshen-Stone/p/book/9780412048418}}

@article{rousseeuw1987,author={Rousseeuw, Peter J.},title={Silhouettes: A Graphical Aid to the Interpretation and Validation of Cluster Analysis},journal={Journal of Computational and Applied Mathematics},volume={20},pages={53--65},year={1987},url={https://doi.org/10.1016/0377-0427(87)90125-7}}

@article{aitchison1982,author={Aitchison, John},title={The Statistical Analysis of Compositional Data},journal={Journal of the Royal Statistical Society: Series B (Methodological)},volume={44},number={2},pages={139--160},year={1982},url={https://doi.org/10.1111/j.2517-6161.1982.tb01195.x}}

@article{lin1991,author={Lin, Jianhua},title={Divergence Measures Based on the {Shannon} Entropy},journal={IEEE Transactions on Information Theory},volume={37},number={1},pages={145--151},year={1991},url={https://doi.org/10.1109/18.61115}}

@article{gneiting2007,author={Gneiting, Tilmann and Raftery, Adrian E.},title={Strictly Proper Scoring Rules, Prediction, and Estimation},journal={Journal of the American Statistical Association},volume={102},number={477},pages={359--378},year={2007},url={https://doi.org/10.1198/016214506000001437}}

@article{pedregosa2011,author={Pedregosa, Fabian and Varoquaux, Ga{\"e}l and Gramfort, Alexandre and Michel, Vincent and Thirion, Bertrand and Grisel, Olivier and Blondel, Mathieu and Prettenhofer, Peter and Weiss, Ron and Dubourg, Vincent and Vanderplas, Jake and Passos, Alexandre and Cournapeau, David and Brucher, Matthieu and Perrot, Matthieu and Duchesnay, {\'E}douard},title={Scikit-learn: Machine Learning in {Python}},journal={Journal of Machine Learning Research},volume={12},pages={2825--2830},year={2011},url={https://jmlr.org/papers/v12/pedregosa11a.html}}

@inproceedings{seabold2010,author={Seabold, Skipper and Perktold, Josef},title={Statsmodels: Econometric and Statistical Modeling with {Python}},booktitle={Proceedings of the 9th Python in Science Conference},pages={92--96},year={2010},url={https://doi.org/10.25080/Majora-92bf1922-011}}

@article{bergmeir2018,author={Bergmeir, Christoph and Hyndman, Rob J. and Koo, Bonsoo},title={A Note on the Validity of Cross-Validation for Evaluating Autoregressive Time Series Prediction},journal={Computational Statistics \& Data Analysis},volume={120},pages={70--83},year={2018},url={https://doi.org/10.1016/j.csda.2017.11.003}}

@article{breiman2001,author={Breiman, Leo},title={Random Forests},journal={Machine Learning},volume={45},pages={5--32},year={2001},url={https://doi.org/10.1023/A:1010933404324}}

@article{roberts2017,author={Roberts, David R. and Bahn, Volker and Ciuti, Simone and Boyce, Mark S. and Elith, Jane and Guillera-Arroita, Gurutzeta and Hauenstein, Severin and Lahoz-Monfort, Jos{\'e} J. and Schr{\"o}der, Boris and Thuiller, Wilfried and Warton, David I. and Wintle, Brendan A. and Hartig, Florian and Dormann, Carsten F.},title={Cross-Validation Strategies for Data with Temporal, Spatial, Hierarchical, or Phylogenetic Structure},journal={Ecography},volume={40},number={8},pages={913--929},year={2017},url={https://doi.org/10.1111/ecog.02881}}

@article{kunsch1989,author={K{\"u}nsch, Hans R.},title={The Jackknife and the Bootstrap for General Stationary Observations},journal={The Annals of Statistics},volume={17},number={3},pages={1217--1241},year={1989},url={https://doi.org/10.1214/aos/1176347265}}

@article{hubert1985,author={Hubert, Lawrence and Arabie, Phipps},title={Comparing Partitions},journal={Journal of Classification},volume={2},pages={193--218},year={1985},url={https://doi.org/10.1007/BF01908075}}

@article{strobl2008,author={Strobl, Carolin and Boulesteix, Anne-Laure and Kneib, Thomas and Augustin, Thomas and Zeileis, Achim},title={Conditional Variable Importance for Random Forests},journal={BMC Bioinformatics},volume={9},pages={307},year={2008},url={https://doi.org/10.1186/1471-2105-9-307}}

@article{cawley2010,author={Cawley, Gavin C. and Talbot, Nicola L. C.},title={On Over-fitting in Model Selection and Subsequent Selection Bias in Performance Evaluation},journal={Journal of Machine Learning Research},volume={11},pages={2079--2107},year={2010},url={https://jmlr.org/papers/v11/cawley10a.html}}

@article{bengio2004,author={Bengio, Yoshua and Grandvalet, Yves},title={No Unbiased Estimator of the Variance of K-Fold Cross-Validation},journal={Journal of Machine Learning Research},volume={5},pages={1089--1105},year={2004},url={https://jmlr.org/papers/v5/grandvalet04a.html}}

@article{hoerl1970,author={Hoerl, Arthur E. and Kennard, Robert W.},title={Ridge Regression: Biased Estimation for Nonorthogonal Problems},journal={Technometrics},volume={12},number={1},pages={55--67},year={1970},url={https://doi.org/10.1080/00401706.1970.10488634}}

@article{stone1977,author={Stone, Charles J.},title={Consistent Nonparametric Regression},journal={The Annals of Statistics},volume={5},number={4},pages={595--620},year={1977},url={https://doi.org/10.1214/aos/1176343886}}

@book{villani2009,author={Villani, C{\'e}dric},title={Optimal Transport: Old and New},publisher={Springer},series={Grundlehren der mathematischen Wissenschaften},volume={338},year={2009},url={https://doi.org/10.1007/978-3-540-71050-9}}
\endgroup
\clearpage
\appendix
\section{Data And Experimental Protocol}
\label{app:protocol}
\subsection{Data Cleaning}
We use the daily Wordle records supplied with 2023 MCM Problem C~\citep{comap2023}. Table~\ref{tab:cleaning} lists the excluded word-level records and the report-count correction.

\begin{table}[!ht]
\caption{Cleaning decisions. The four word-level exclusions retain their dates in the count series. Other admissible rounded histograms are normalized to 100\%.}
\label{tab:cleaning}\centering\small
\begin{tabular*}{\linewidth}{@{\extracolsep{\fill}}lll@{}}
\toprule
Record & Issue & Action\\
\midrule
2022-03-27, NYMPH & Outcome sum 126\% & Exclude word-level row\\
2022-04-29, TASH & Four-letter string & Exclude word-level row\\
2022-11-26, CLEN & Four-letter string & Exclude word-level row\\
2022-12-16, RPROBE & Six-letter string & Exclude word-level row\\
FAVOR / NA\"IVE & Whitespace / diacritic & Canonicalize to five letters\\
2022-11-30 & Total 2,569; hard mode 2,405 & Repair total to 25,546.5\\
\bottomrule
\end{tabular*}
\end{table}

Malformed solution words are excluded without replacement. The November 30 total is imputed from the six preceding daily counts. Outcome histograms are retained when their total differs from 100\% by at most 3.5 pp, the largest discrepancy permitted by rounding seven shares to integer percentages.

\begin{table}[!ht]
\caption{Exact temporal partitions. Counts are fitting/validation/test observations before refitting on fitting plus validation. All dates are in 2022.}
\label{tab:splits}\centering\small
\begin{tabular*}{\linewidth}{@{\extracolsep{\fill}}llrrr@{}}
\toprule
Split & Fit/Val/Test & Fitting Dates & Validation Dates & Test Dates\\
\midrule
Primary & 198/50/107 & 01-07--07-25 & 07-26--09-13 & 09-14--12-31\\
F1 & 124/31/50 & 01-07--05-12 & 05-13--06-12 & 06-13--08-01\\
F2 & 164/41/50 & 01-07--06-21 & 06-22--08-01 & 08-02--09-20\\
F3 & 204/51/50 & 01-07--07-31 & 08-01--09-20 & 09-21--11-09\\
F4 & 244/61/50 & 01-07--09-09 & 09-10--11-09 & 11-10--12-31\\
\bottomrule
\end{tabular*}

\end{table}

\subsection{Nested Selection And Comparator Settings}
Table~\ref{tab:splits} gives every boundary. In each split, feature standardization for ridge and kNN is fitted only on the current fitting partition during selection, then refitted on the full training prefix. The last 20\% of the prefix is the validation segment (integer fitting size $\lfloor0.8n\rfloor$). Selection minimizes validation MAE after clipping and normalization; ties use the first candidate in the recorded grid. The outer test never determines a parameter, transform, or cluster center. Later expanding-window training prefixes may include earlier folds' test words, which have become historical observations by that later origin.

All candidates use the same three lexical inputs. Ridge selects $\alpha\in\{0.1,1,10,100\}$ with an unpenalized intercept. kNN uses Euclidean distance in standardized feature space, inverse-distance weights, and $k\in\{10,25,50\}$; exact input matches follow the library's zero-distance weighting convention. Random forests use 300 bootstrap trees, all three candidate features at each split, maximum depth in $\{2,4\}$, and minimum leaf size in $\{5,10\}$. They regress the seven outputs jointly. XGBoost selects maximum depth in $\{1,2,3\}$ and rounds in $\{100,300\}$; its six independent regressors and constant raw first-bin share follow Section~\ref{sec:xgb}. Full row/column sampling, $\eta=0.03$, $\lambda=10$, $\gamma=0$, no $\ell_1$ penalty, histogram tree construction, and minimum child weight 5 remain fixed. The mean comparator averages training histograms. No comparator receives additional word information.

\begin{table}[!ht]
\caption{Selected parameters in each temporal split. Ridge entries are $\alpha$; kNN entries are $k$; forest entries are depth/minimum leaf size; XGBoost entries are depth/rounds.}
\label{tab:selected}\centering\small
\begin{tabular*}{\linewidth}{@{\extracolsep{\fill}}lrrrrr@{}}
\toprule
Model & Primary & F1 & F2 & F3 & F4\\
\midrule
Ridge & 0.1 & 100 & 10 & 0.1 & 0.1\\
kNN & 50 & 50 & 50 & 50 & 50\\
Random Forest & 2/5 & 2/10 & 2/10 & 2/5 & 2/10\\
XGBoost & 1/300 & 1/100 & 2/100 & 1/300 & 2/100\\
\bottomrule
\end{tabular*}

\end{table}

\begin{table}[!ht]
\caption{Complete primary XGBoost validation grid, MAE (pp). The selected entry is bold. These are validation errors, not test scores.}
\label{tab:grid}\centering\small
\begin{tabular*}{\linewidth}{@{\extracolsep{\fill}}lrr@{}}
\toprule
\textbf{Depth} & \textbf{100 Rounds} & \textbf{300 Rounds}\\
\midrule
1 & 3.540 & \textbf{3.511}\\
2 & 3.514 & 3.582\\
3 & 3.607 & 3.636\\
\bottomrule
\end{tabular*}

\end{table}

Table~\ref{tab:selected} reports the selected configurations, and Table~\ref{tab:grid} gives the primary XGBoost validation scores.

\subsection{Implementation}
The experiments use Python 3.14.6, scikit-learn 1.9.1~\citep{pedregosa2011}, XGBoost 3.4.1, statsmodels 0.15.0~\citep{seabold2010}, and wordfreq 3.1.1. Tree models use one worker and random seed 42 unless otherwise specified.

\section{Formalization And Evaluation Metrics}
\label{app:metrics}
\subsection{Estimand And Boosting Updates}
The observable histogram concerns reporting players. An explicit population target at date $t$ is
\begin{equation}
 p^*_{t,b}(w)=100\,\Pr(G_t=b\mid W_t=w,S_t=1),\qquad b\in\mathcal B,
 \label{eq:estimand}
\end{equation}
where $G_t$ is a player's outcome, $W_t$ is the solution word, and $S_t$ indicates that the player publicly reports a result. The workbook supplies rounded aggregate estimates of these percentages, not player-level draws. Lexical regression pools across dates to approximate the histogram from $x(w)$, and does not separately identify a word effect from temporal changes in the reporting population. Equal word weights target an average word-level error, not an average error over individual reports.

For the half-squared loss in Equation~\eqref{eq:xgb-objective}, a candidate tree partition has leaf updates
\begin{equation}
 g_i=z_b^{(m-1)}(x_i)-p_{i,b},\qquad h_i=1,\qquad
 a_\ell^*=-\frac{\sum_{i\in\mathcal I_\ell}g_i}{\sum_{i\in\mathcal I_\ell}h_i+\lambda},
 \label{eq:leaf-update}
\end{equation}
where $\mathcal I_\ell$ contains the fitting observations routed to leaf $\ell$, $g_i$ and $h_i$ are the loss gradient and Hessian before the added tree, and $a_\ell^*$ is the regularized pre-shrinkage leaf weight. The implementation adds $\eta a_\ell^*$ to the raw score. Equation~\eqref{eq:leaf-update} is the exact leaf optimum for a fixed partition; split search is approximate~\citep{chen2016,friedman2001}.

\subsection{Distribution And Classification Metrics}
The primary average bin error and the ordered distribution error are
\begin{equation}
 \begin{aligned}
 \operatorname{MAE}_{\rm pp}&=\frac{1}{7|\mathcal E|}\sum_{i\in\mathcal E}\sum_{b\in\mathcal B}|p_{i,b}-\widehat p_{i,b}|,\\
 W_1(u,v)&=\sum_{j=1}^{6}\left|\sum_{b=1}^{j}(u_b-v_b)\right|,
 \end{aligned}
 \label{eq:distribution-errors}
\end{equation}
where $\mathcal E$ is the evaluated test index set, $u=\mathbf p_i/100$, $v=\widehat{\mathbf p}_i/100$, and failure occupies ordered position 7 for $W_1$. MAE is measured in pp. We report the mean $W_1$ over test words on the unit-spaced capped outcome scale~\citep{villani2009}.

The remaining probability-vector discrepancies are
\begin{equation}
 \begin{aligned}
 \operatorname{TV}(u,v)&=\tfrac12\sum_{b\in\mathcal B}|u_b-v_b|,\\
 \operatorname{JSD}(u,v)&=\tfrac12\operatorname{KL}_2(u\|m)+\tfrac12\operatorname{KL}_2(v\|m),\qquad m=(u+v)/2,
 \end{aligned}
 \label{eq:distribution-metrics}
\end{equation}
where $\operatorname{KL}_2(u\|m)=\sum_bu_b\log_2(u_b/m_b)$ and zero-numerator terms are zero. TV and base-2 JSD lie in $[0,1]$, with JSD measured in bits. Under equal weights across the seven percentage bins, mean TV equals $7\operatorname{MAE}_{\rm pp}/200$. The main comparison instead reports the complementary $W_1$ metric.

The marginal tolerance score is
\begin{equation}
 A_3=\frac{100}{6|\mathcal E|}\sum_{i\in\mathcal E}\sum_{b\in\mathcal B_+}\mathbf1\{|p_{i,b}-\widehat p_{i,b}|\leq3\},
 \label{eq:tolerance}
\end{equation}
where $\mathcal B_+=\{2,3,4,5,6,\fail\}$ excludes the raw constant-prediction first bin. All predictions are evaluated after normalization. For the primary mean, ridge, and XGBoost models, $A_3$ is 45.95\%, 51.71\%, and 49.07\%, respectively. It averages individual bin events and is not a seven-bin confidence level.

Classification accuracy is the fraction of correctly predicted reference tiers. Balanced accuracy averages the $K$ class recalls. Macro-F1 averages $2\mathrm{TP}_k/(2\mathrm{TP}_k+\mathrm{FP}_k+\mathrm{FN}_k)$ over classes, using zero for an undefined class score. The reference labels come from the fixed training centroids. The majority comparator predicts the most common \emph{training} tier for all test words.

\subsection{Paired Uncertainty And Dependence}
For a fitted model $A$ and comparator $B$, define per-word loss differences and their bootstrap means by
\begin{equation}
 d_i=\tfrac17\sum_{b\in\mathcal B}\bigl(|p_{i,b}-\widehat p^A_{i,b}|-|p_{i,b}-\widehat p^B_{i,b}|\bigr),\qquad
 \overline d^{*(r)}=\tfrac1n\sum_{j=1}^{n}d_{I_j^{(r)}},
 \label{eq:paired-bootstrap}
\end{equation}
where $n=107$ and $I^{(r)}$ is a length-$n$ index sequence formed by concatenating uniformly sampled circular blocks of length $L$ and truncating to $n$. We use 2,000 resamples and seed 42. The empirical 2.5th and 97.5th percentiles give the interval for the period-specific paired mean. Circular blocks preserve short-range dependence under approximate local stability~\citep{kunsch1989}.

\begin{table}[!ht]
\caption{Sensitivity of paired MAE-difference intervals to block length. Each model is compared with the training mean; units are pp. The mean's own difference is identically zero.}
\label{tab:blocks}\centering\small
\begin{tabular*}{\linewidth}{@{\extracolsep{\fill}}lrrr@{}}
\toprule
Model & Block 1 & Block 7 & Block 14\\
\midrule
Training Mean & $[+0.000,\,+0.000]$ & $[+0.000,\,+0.000]$ & $[+0.000,\,+0.000]$\\
Ridge & $[-0.502,\,-0.067]$ & $[-0.490,\,-0.062]$ & $[-0.462,\,-0.102]$\\
kNN & $[-0.174,\,+0.550]$ & $[-0.120,\,+0.490]$ & $[-0.065,\,+0.420]$\\
Random Forest & $[-0.539,\,-0.119]$ & $[-0.522,\,-0.134]$ & $[-0.480,\,-0.184]$\\
XGBoost & $[-0.445,\,-0.044]$ & $[-0.424,\,-0.036]$ & $[-0.395,\,-0.096]$\\
\bottomrule
\end{tabular*}

\end{table}

Table~\ref{tab:blocks} checks $L\in\{1,7,14\}$. The intervals condition on fitted models and the observed test period. The temporal folds share training observations~\citep{bengio2004}.

The primary improvement over the training mean is not driven by a single test word. Holding fitted predictions fixed, deleting any one word's loss leaves the XGBoost-minus-mean difference between $-0.269$ and $-0.220$ pp. XGBoost improves per-word MAE on 60 of 107 words. Predictions remain fixed during this influence analysis.

\section{Detailed Word-Level Results}
\label{app:errors}
\subsection{Subsets And Signed Errors}
\begin{table}[!ht]
\caption{Primary holdout slices. Model entries are MAE (pp). XGB X-Bias is the mean predicted-minus-observed failure share (pp). Subsets overlap; they are descriptive views, not separate test sets.}
\label{tab:slices}\centering\small
\begin{tabular*}{\linewidth}{@{\extracolsep{\fill}}lrrrrrrr@{}}
\toprule
Subset & $n$ & Mean & Ridge & kNN & RF & XGB & XGB X-Bias (pp)\\
\midrule
All & 107 & 4.082 & 3.799 & 4.262 & 3.750 & 3.839 & +0.05\\
Rare & 55 & 4.135 & 3.986 & 4.562 & 3.979 & 4.073 & -0.22\\
Frequent & 52 & 4.027 & 3.601 & 3.945 & 3.508 & 3.592 & +0.33\\
Distinct & 81 & 3.946 & 3.696 & 4.216 & 3.634 & 3.716 & +0.81\\
Repeated & 26 & 4.508 & 4.119 & 4.406 & 4.113 & 4.222 & -2.33\\
High failure & 8 & 8.333 & 7.501 & 7.982 & 7.668 & 7.735 & -15.25\\
\bottomrule
\end{tabular*}

\end{table}

We define rare and frequent words using the median log-frequency of the 248 training words, $-5.4202$; ties belong to the frequent group. Repeated denotes $\NRE>0$ and distinct denotes $\NRE=0$. High failure uses the training 90th percentile, 6.0606 pp, with the test criterion $p_{i,\fail}\geq6.0606$. The high-failure subset uses observed test outcomes and is intended for error analysis.

The direction persists beyond the primary cutoff. Recomputing the 90th-percentile threshold from each temporal fold's training prefix gives 7, 8, 4, and 2 high-failure test words in F1--F4. Their XGBoost failure-bin biases are $-6.87$, $-11.07$, $-12.96$, and $-9.87$ pp, respectively. In every fold, tail MAE also exceeds the remaining words' MAE. The tail subsets are small, but their failure-rate biases have the same direction across all four windows.

\begin{table}[!ht]
\caption{Per-bin MAE (pp) and XGBoost signed bias (pp) on all 107 test words. The near-zero aggregate failure bias can coexist with severe underprediction of rare high-failure words.}
\label{tab:perbin}\centering\small
\begin{tabular*}{\linewidth}{@{\extracolsep{\fill}}lrrrrrr@{}}
\toprule
Outcome & Mean & Ridge & kNN & RF & XGB & XGB Bias\\
\midrule
1 & 0.649 & 0.594 & 0.513 & 0.580 & 0.649 & +0.117\\
2 & 3.027 & 2.766 & 3.019 & 2.699 & 2.776 & -0.585\\
3 & 7.087 & 6.254 & 7.291 & 6.286 & 6.340 & -2.504\\
4 & 4.577 & 4.486 & 4.933 & 4.472 & 4.658 & -1.214\\
5 & 5.302 & 5.077 & 5.525 & 4.836 & 4.896 & +1.914\\
6 & 5.339 & 4.881 & 5.887 & 4.868 & 5.018 & +2.226\\
X & 2.595 & 2.533 & 2.664 & 2.510 & 2.535 & +0.047\\
\bottomrule
\end{tabular*}

\end{table}

\begin{figure}[!ht]
\centering
\begin{subfigure}[t]{0.32\linewidth}
\centering
\includegraphics[width=\linewidth]{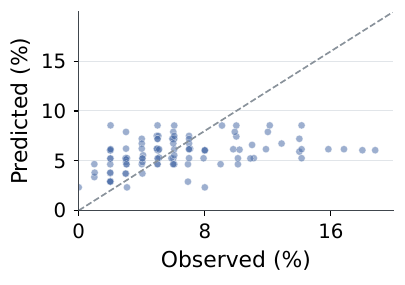}
\caption{2 attempts}\label{fig:scatter:a}
\end{subfigure}
\hfill
\begin{subfigure}[t]{0.32\linewidth}
\centering
\includegraphics[width=\linewidth]{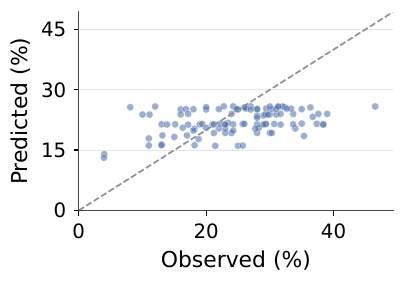}
\caption{3 attempts}\label{fig:scatter:b}
\end{subfigure}
\hfill
\begin{subfigure}[t]{0.32\linewidth}
\centering
\includegraphics[width=\linewidth]{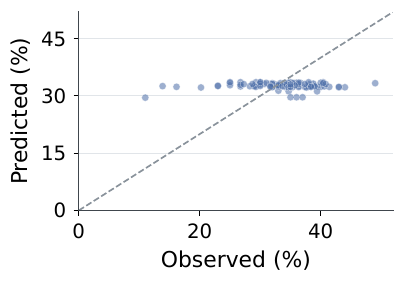}
\caption{4 attempts}\label{fig:scatter:c}
\end{subfigure}
\par
\begin{subfigure}[t]{0.32\linewidth}
\centering
\includegraphics[width=\linewidth]{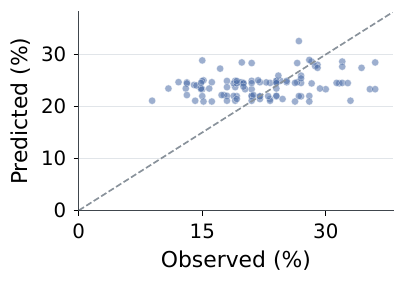}
\caption{5 attempts}\label{fig:scatter:d}
\end{subfigure}
\hfill
\begin{subfigure}[t]{0.32\linewidth}
\centering
\includegraphics[width=\linewidth]{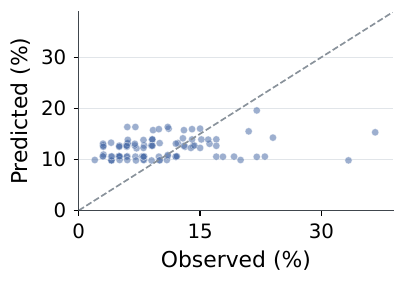}
\caption{6 attempts}\label{fig:scatter:e}
\end{subfigure}
\hfill
\begin{subfigure}[t]{0.32\linewidth}
\centering
\includegraphics[width=\linewidth]{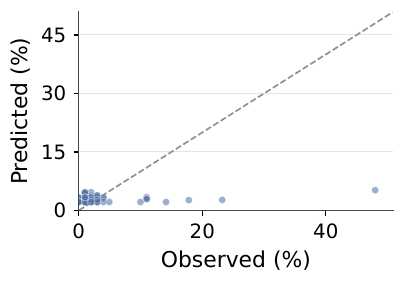}
\caption{Failure (X)}\label{fig:scatter:f}
\end{subfigure}
\caption{Observed and predicted outcome shares for held-out words. Dashed lines indicate exact agreement.}
\label{fig:scatter}
\end{figure}

Table~\ref{tab:perbin} explains why a single average bias is insufficient. In the full test set, XGBoost failure-bin bias is approximately $+0.05$ pp, yet it is $-15.25$ pp within the high-failure subset. Positive errors on common low-failure words offset negative errors on difficult tails. The observed/predicted plots in Figure~\ref{fig:scatter} reveal this shrinkage toward typical histograms across outcomes.

\subsection{Feature Reliance And Capacity}
\begin{figure}[!ht]
\centering
\begin{subfigure}[t]{0.32\linewidth}
\centering
\includegraphics[width=\linewidth]{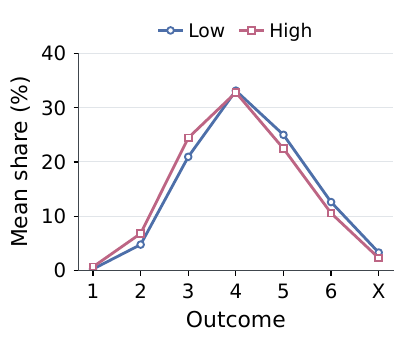}
\caption{FREQ}\label{fig:profiles:a}
\end{subfigure}
\hfill
\begin{subfigure}[t]{0.32\linewidth}
\centering
\includegraphics[width=\linewidth]{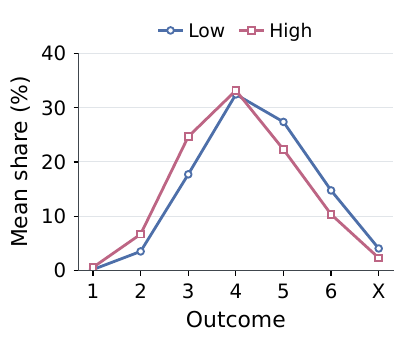}
\caption{WIE}\label{fig:profiles:b}
\end{subfigure}
\hfill
\begin{subfigure}[t]{0.32\linewidth}
\centering
\includegraphics[width=\linewidth]{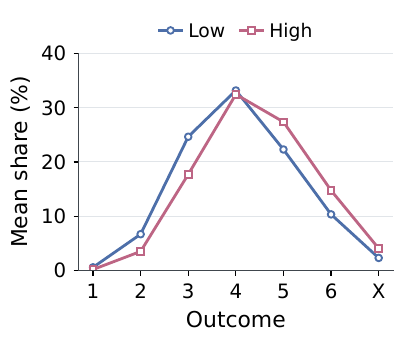}
\caption{NRE}\label{fig:profiles:c}
\end{subfigure}
\caption{Average outcome distributions for high and low attribute values. Low/high group sizes are 177/178 for FREQ, 100/255 for WIE, and 255/100 for NRE.}
\label{fig:profiles}
\end{figure}

Figure~\ref{fig:profiles} splits FREQ at its full-dataset median; ties belong to high. High WIE corresponds to five distinct letters, and high NRE to any repeated letter. Low/high sizes are 177/178, 100/255, and 255/100. These descriptive groups are not used for predictive preprocessing or model selection. FREQ correlates 0.36 with two-guess success and $-0.23$ with failure; NRE correlates $-0.46$ with two-guess success and 0.25 with failure. These exploratory associations carry no causal claim.

The ablation retains the primary selected XGBoost depth and number of rounds for all feature subsets. It evaluates representation changes under common capacity, not separately optimized feature-set performance. Ten random 70/30 splits with this configuration give MAE $3.429\pm0.124$ pp. Because the configuration was selected chronologically, these random splits may reuse selection observations in their test sets; the nested temporal tests provide the primary assessment of generalization.

The 355 words exhibit four letter-count patterns: $1+1+1+1+1$ (255 words), $2+1+1+1$ (94), $2+2+1$ (4), and $3+1+1$ (2). On these observed patterns, WIE and NRE determine one another exactly. Their different numerical orderings can change tree partitions, but neither identifies letter identities or positions. The repeated-letter category therefore retains substantial outcome variation.

\begin{figure}[!ht]
\centering
\begin{subfigure}[t]{.32\linewidth}\centering
\includegraphics[width=\linewidth]{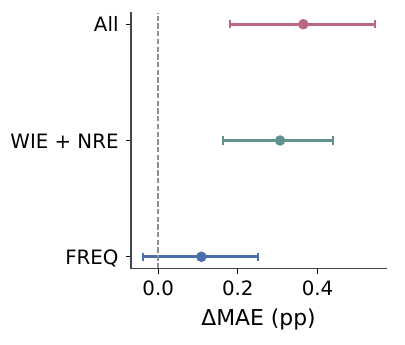}\caption{Grouped permutation}
\end{subfigure}\hfill
\begin{subfigure}[t]{.32\linewidth}\centering
\includegraphics[width=\linewidth]{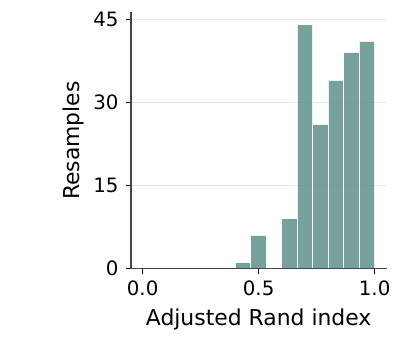}\caption{Partition agreement}
\end{subfigure}\hfill
\begin{subfigure}[t]{.32\linewidth}\centering
\includegraphics[width=\linewidth]{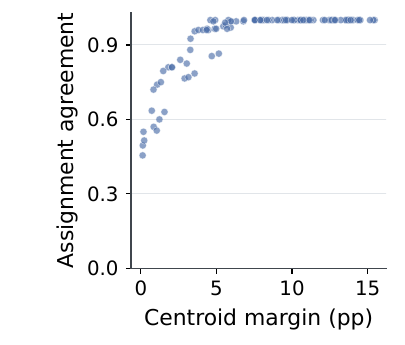}\caption{Boundary sensitivity}
\end{subfigure}
\caption{Model reliance and cluster sensitivity. Bars in (a) span the 2.5th--97.5th percentiles of 200 permutation effects, not sampling confidence intervals. Panel (b) compares resampled partitions with the reference using ARI; (c) relates per-word assignment agreement to its distance margin.}
\label{fig:extended-diagnostics}
\end{figure}

Permutation diagnostics hold all fitted parameters fixed. Each repetition permutes test rows within a feature group, applies that row permutation jointly to its columns, and recomputes MAE. Jointly permuting WIE and NRE avoids breaking their internal relationship, but still breaks association with other groups. The empirical permutation range describes variation across shuffles, not uncertainty across new datasets. Correlated-feature importance requires particular caution~\citep{strobl2008}; neither this display nor tree impurity importance identifies independent causal contributions.

\subsection{Largest-Error Cases}
\label{app:cases}
\begin{table}[!ht]
\caption{Ten largest primary XGBoost histogram errors, ordered by MAE (pp). Failure shares are percentages; $d$ is the capped score from Section~\ref{sec:tiers}. Cases are selected after evaluation for diagnosis.}
\label{tab:cases}\centering\small
\begin{tabular*}{\linewidth}{@{\extracolsep{\fill}}llrrrrr@{}}
\toprule
Date & Word & MAE & Observed X & Predicted X & Observed $d$ & Predicted $d$\\
\midrule
2022-09-16 & PARER & 12.921 & 48.00 & 5.20 & 5.99 & 4.63\\
2022-12-26 & JUDGE & 11.661 & 14.14 & 2.16 & 5.23 & 4.03\\
2022-10-23 & MUMMY & 10.432 & 17.82 & 2.65 & 5.48 & 4.50\\
2022-12-05 & WOKEN & 7.703 & 4.00 & 2.17 & 4.80 & 4.13\\
2022-09-17 & CHUTE & 7.285 & 0.00 & 2.74 & 3.59 & 4.25\\
2022-11-28 & TEPID & 7.015 & 0.00 & 3.31 & 3.61 & 4.27\\
2022-10-15 & CATCH & 7.004 & 23.23 & 2.70 & 4.85 & 4.25\\
2022-11-04 & PHOTO & 6.966 & 0.00 & 2.71 & 3.77 & 4.25\\
2022-10-06 & SLOTH & 6.891 & 0.00 & 2.95 & 3.58 & 4.25\\
2022-11-09 & RAINY & 6.799 & 0.99 & 2.19 & 3.48 & 4.13\\
\bottomrule
\end{tabular*}

\end{table}

Table~\ref{tab:cases} includes both unusually difficult words (PARER, JUDGE, MUMMY) and overpredicted difficulty on CHUTE, TEPID, and other words. PARER has 48\% observed failure versus a 5.20\% prediction, illustrating how a lexically compact predictor misses a distinctive outcome profile. The ranking includes every word in the ten largest errors. The observed features do not identify whether these errors arise from morphology, player strategies, or lexical neighborhoods.

The \EERIE{} plot uses the primary 248-word model. Refitting XGBoost on all 355 retained words gives $(0.46,2.92,11.11,29.64,33.09,19.34,3.44)$\%. This is an illustrative case prediction without an observed target. Figure~\ref{fig:eerie}c additionally displays the concentration statistic $s_i=\sum_{b\in\{3,4,5,6,\fail\}}p_{i,b}$. Exactly 299 of 355 words satisfy $s_i>90$, while three equal 90. The plotted empirical CDF $\widehat F(s)=355^{-1}\sum_i\mathbf1\{s_i\leq s\}$ gives the strict fraction $1-\widehat F(90)=84.225\%$. This is a fraction of word-level distributions, not a fraction of all Wordle players.

\section{Difficulty Prediction And Stability}
\label{app:stability}
\label{app:routes}
\subsection{Direct And Indirect Prediction}
The primary lexical tree predicts a tier directly from the three attributes. An indirect route first predicts a histogram with a fitted regressor, then assigns it to its nearest fixed training centroid. Table~\ref{tab:routes} evaluates both against the same 107 reference labels, without fitting additional models. The regressors' predictions are those already evaluated in Table~\ref{tab:prediction}.

\begin{table}[!ht]
\caption{Primary difficulty routes using the same three training centroids. Accuracy, balanced accuracy, and macro-F1 are percentages. Reference counts are 53 easy, 36 moderate, and 18 difficult; the final column counts predicted labels.}
\label{tab:routes}\centering\small
\begin{tabular*}{\linewidth}{@{\extracolsep{\fill}}lrrrr@{}}
\toprule
Route & Accuracy & Bal. Acc. & Macro-F1 & Predicted E/M/D\\
\midrule
Lexical Tree & \textbf{54.2} & \textbf{51.3} & \textbf{50.4} & 50/31/26\\
Mean Histogram & 33.6 & 33.3 & 16.8 & 0/107/0\\
Ridge Histogram & 39.3 & 37.8 & 31.9 & 16/86/5\\
kNN Histogram & 38.3 & 38.7 & 35.4 & 24/62/21\\
Random Forest Histogram & 38.3 & 40.8 & 34.5 & 10/87/10\\
XGBoost Histogram & 40.2 & 39.0 & 32.4 & 15/90/2\\
\bottomrule
\end{tabular*}

\end{table}

The tree predicts 50 easy, 31 moderate, and 26 difficult words. The XGBoost route instead predicts 15, 90, and 2, respectively. Ridge and random forests show similar concentration in the moderate tier. Even the forest, which has the lowest primary histogram MAE, attains only 38.3\% tier accuracy. Average histogram fidelity and the discontinuous nearest-centroid decision therefore evaluate different behavior. The direct classifier gives higher tier accuracy despite using the same lexical inputs.

\subsection{Resampling And Decision Margins}
For each of 200 resamples, we draw 248 training profiles with replacement, fit $K=3$ using 50 k-means++ restarts, and assign the fixed 107 test profiles to those centers. The resampling generator uses seed 42 and K-means uses the resample index as its seed. Centers are ordered by capped score, as in the reference fit. This procedure measures sensitivity of the tier assignments to the training sample.

Let $c_i^0$ be a reference assignment and $c_i^{(r)}$ its assignment under resample $r$. Per-word stability and the distance margin are
\begin{equation}
 a_i=\frac1R\sum_{r=1}^{R}\mathbf1\{c_i^{(r)}=c_i^0\},\qquad
 \delta_i=D_{i,(2)}-D_{i,(1)},
 \label{eq:stability}
\end{equation}
where $R=200$ and $D_{i,(1)}\leq D_{i,(2)}$ are the smallest two Euclidean distances from test profile $\mathbf p_i$ to the reference training centroids. The margin has pp units in histogram space. Smaller margins indicate weaker geometric separation, not statistical confidence. The mean $a_i$ is 0.930; 16 words have $a_i<0.8$ (Figure~\ref{fig:extended-diagnostics}c).

This margin also gives a sufficient condition for the indirect route to preserve the reference tier:
\begin{equation}
 e_i<\frac{\delta_i}{2}
 \quad\Longrightarrow\quad
 \underset{k\in\{1,2,3\}}{\arg\min}\,
 \|\widehat{\mathbf p}_i-\mu_k\|_2=c_i^0,
 \label{eq:tier-margin}
\end{equation}
where $e_i=\|\widehat{\mathbf p}_i-\mathbf p_i\|_2$ is the histogram error in pp and $\delta_i$ is the reference distance gap from Equation~\eqref{eq:stability}. By the triangle inequality, moving a profile by $e_i$ increases its distance to the winning centroid by at most $e_i$ and decreases its distance to a competitor by at most $e_i$. A gap greater than $2e_i$ therefore preserves the winner. The condition is sufficient, not necessary: it holds for 4 of 107 XGBoost predictions, whereas 43 receive the correct tier. Because $e_i$ and $\delta_i$ require the observed histogram, this is an explanation of evaluated decisions, not a confidence guarantee available for an unseen word.

Adjusted Rand index (ARI) is invariant to label permutations and corrects co-membership agreement for chance~\citep{hubert1985}: its contingency-table expression is
\begin{equation}
 \operatorname{ARI}=\frac{\sum_{a,b}\binom{n_{ab}}2-AB/\binom n2}
 {\tfrac12(A+B)-AB/\binom n2},\qquad
 A=\sum_a\binom{n_{a+}}2,\quad B=\sum_b\binom{n_{+b}}2,
 \label{eq:ari}
\end{equation}
where $n_{ab}$ counts test words in reference cluster $a$ and resampled cluster $b$, $n_{a+}$ and $n_{+b}$ are the corresponding margins, and $n=107$. The denominator is nonzero for the evaluated partitions. ARI is 1 for identical partitions and can be negative when agreement is below its chance reference. This correction differs from the directly interpretable ordered-label agreement $a_i$.

\subsection{Tier Granularity}
\begin{table}[!ht]
\caption{Tier-count sensitivity with the same depth-3 lexical tree and minimum leaf size 10. Metrics are percentages. Each $K$ defines a different target; test scores are not used to choose $K$.}
\label{tab:tiersensitivity}\centering\small
\begin{tabular*}{\linewidth}{@{\extracolsep{\fill}}rrrrr@{}}
\toprule
$K$ & Accuracy & Balanced Accuracy & Macro-F1 & Majority Accuracy\\
\midrule
2 & 57.0\% & 59.7\% & 55.9\% & 69.2\%\\
3 & 54.2\% & 51.3\% & 50.4\% & 33.6\%\\
4 & 38.3\% & 28.2\% & 24.1\% & 43.0\%\\
5 & 43.9\% & 25.6\% & 20.9\% & 43.0\%\\
\bottomrule
\end{tabular*}

\end{table}

The primary ordered centroid scores are 3.826, 4.237, 4.706 for easy, moderate, and difficult. Test reference counts are 53, 36, and 18. The tree has 53.6\% training accuracy and 54.2\% test accuracy. It splits on FREQ and WIE, with zero NRE impurity importance; this is consistent with the redundancy between WIE and NRE. Table~\ref{tab:tiersensitivity} shows how class balance and accuracy vary with $K$. At $K=2$, the majority baseline exceeds tree accuracy. We retain three tiers for interpretability.

\section{Report Forecasts And Interval Diagnostics}
\label{app:forecasts}
\subsection{Horizon-Matched Protocol And Scores}
We use 18 weekly origins from July 4 through October 31, 2022. At each origin, the models are fitted to counts observed through that date and issue forecasts at all five horizons before later counts are observed. The full model uses all history from January 7; recent models use the 90 preceding observations including the origin. All 54 ARIMA fits converged. We do not tune the window, order, or drift using these target observations.

Count-series point errors are
\begin{equation}
 \operatorname{MAE}_r=\frac1{|\mathcal D|}\sum_{t\in\mathcal D}|r_t-\widehat r_t|,\qquad
 \operatorname{MAPE}_r=\frac{100}{|\mathcal D|}\sum_{t\in\mathcal D}\left|\frac{r_t-\widehat r_t}{r_t}\right|,
 \label{eq:count-metrics}
\end{equation}
where $\mathcal D$ is the set of targets under the stated evaluation protocol and $r_t>0$ is the repaired observed total. The one-step sequential protocol has 61 targets, whereas each matched horizon has 18. Width averages $u-l$, coverage averages $\mathbf1\{l\leq r_t\leq u\}$, and the interval score (lower is better) is
\begin{equation}
 \operatorname{IS}_{\alpha}(l,u;y)=(u-l)+\frac{2}{\alpha}(l-y)\mathbf1\{y<l\}
 +\frac{2}{\alpha}(y-u)\mathbf1\{y>u\},
 \label{eq:interval-score}
\end{equation}
where $[l,u]$ is a central $(1-\alpha)$ prediction interval, $y$ is the realized count, and $\alpha=0.05$~\citep{gneiting2007}. A smaller score rewards narrower intervals while penalizing misses.

\begin{table}[!ht]
\caption{Matched-horizon report forecasts, 18 origins per horizon. Errors, widths, and interval scores use thousands of reports (k). Negative is the percentage of intervals with a negative lower endpoint. Persistence supplies only point forecasts.}
\label{tab:horizons}\centering\small
\begin{tabular*}{\linewidth}{@{\extracolsep{\fill}}rlrrrrr@{}}
\toprule
$h$ & Model & MAE (k) & Coverage (\%) & Width (k) & IS$_{.05}$ (k) & Negative (\%)\\
\midrule
1 & Persistence & 1.52 & -- & -- & -- & --\\
1 & ARIMA Full & 1.44 & 100.0 & 57.2 & 57.2 & 0.0\\
1 & ARIMA 90d & 1.46 & 100.0 & 13.5 & 13.5 & 0.0\\
1 & Log ARIMA 90d & 1.29 & 100.0 & 8.9 & 8.9 & 0.0\\
7 & Persistence & 2.41 & -- & -- & -- & --\\
7 & ARIMA Full & 2.26 & 100.0 & 116.3 & 116.3 & 100.0\\
7 & ARIMA 90d & 2.33 & 100.0 & 30.6 & 30.6 & 0.0\\
7 & Log ARIMA 90d & 2.13 & 94.4 & 12.5 & 13.4 & 0.0\\
14 & Persistence & 2.51 & -- & -- & -- & --\\
14 & ARIMA Full & 2.36 & 100.0 & 159.7 & 159.7 & 100.0\\
14 & ARIMA 90d & 2.45 & 100.0 & 42.6 & 42.6 & 11.1\\
14 & Log ARIMA 90d & 2.44 & 100.0 & 15.8 & 15.8 & 0.0\\
30 & Persistence & 3.40 & -- & -- & -- & --\\
30 & ARIMA Full & 3.18 & 100.0 & 229.8 & 229.8 & 100.0\\
30 & ARIMA 90d & 3.28 & 100.0 & 61.9 & 61.9 & 22.2\\
30 & Log ARIMA 90d & 3.57 & 100.0 & 21.6 & 21.6 & 0.0\\
60 & Persistence & 6.32 & -- & -- & -- & --\\
60 & ARIMA Full & 6.13 & 100.0 & 322.6 & 322.6 & 100.0\\
60 & ARIMA 90d & 6.23 & 100.0 & 87.2 & 87.2 & 66.7\\
60 & Log ARIMA 90d & 7.04 & 100.0 & 29.9 & 29.9 & 0.0\\
\bottomrule
\end{tabular*}

\end{table}

For log-count ARIMA, the forecast mean and interval on the count scale are
\begin{equation}
 \begin{aligned}
 \widehat r_{T+h\mid T}^{\log}&=\exp(m_{T,h}^{\log}+v_{T,h}^{\log}/2),\\
 I_{T,h}^{\log}&=\left[\exp(m_{T,h}^{\log}-z_{.975}\sqrt{v_{T,h}^{\log}}),\;
 \exp(m_{T,h}^{\log}+z_{.975}\sqrt{v_{T,h}^{\log}})\right],
 \end{aligned}
 \label{eq:log-forecast}
\end{equation}
where $m^{\log}$ and $v^{\log}$ are the fitted conditional mean and variance of the Gaussian log-count forecast. The variance correction converts its median to its mean~\citep{hyndman2021}. These intervals respect positivity, but do not incorporate parameter uncertainty or establish a correct count-generating model. Overlapping targets and one-year nonstationarity limit calibration inference from the 18 origins.

\subsection{Repair Sensitivity}
The matched-origin fits precede November 30, so this repair does not affect their fitting inputs. One 30-day target equals that date. Table~\ref{tab:rawcounts} evaluates the same fixed predictions against the repaired and original totals, making the impact of the chosen target explicit. Results at the other four horizons are unchanged.

\begin{table}[!ht]
\caption{Thirty-day MAE under the two November 30 target choices, thousands of reports. Each column averages all 18 origins; only one target differs. Prediction models are unchanged.}
\label{tab:rawcounts}\centering\small
\begin{tabular*}{\linewidth}{@{\extracolsep{\fill}}lrr@{}}
\toprule
Model & Repaired Target MAE (k) & Original Target MAE (k)\\
\midrule
Persistence & 3.399 & 4.676\\
ARIMA Full & 3.183 & 4.460\\
ARIMA 90d & 3.284 & 4.561\\
Log ARIMA 90d & 3.572 & 4.848\\
\bottomrule
\end{tabular*}

\end{table}

\subsection{Residual And Window Diagnostics}
\label{app:diagnostics}
The reporting trajectory contains a pronounced adoption peak followed by decline. The full-history level ADF $p$-value is 0.001883 and the difference value is 0.000576. Rejection of a unit-root null does not establish a constant mean and variance through these phases. The fitted moving-average coefficient is $-0.27934$; the lag-10 Ljung--Box residual $p$-value is $1.12\times10^{-7}$, inconsistent with residual white noise.

Coefficient sensitivity fixes $\theta\in\{-0.8,-0.6,\ldots,0.8\}$ and re-estimates innovation variance under the same zero-drift sign convention. The July--December refit gives a March 1 mean of 20,527 and interval $[-338,41{,}392]$, compared with the full-history interval $[-109{,}546,150{,}588]$. The matched-horizon tests in Table~\ref{tab:horizons} evaluate this window sensitivity against observed counts.

\end{document}